\documentclass[12pt]{article}
\pdfoutput=1
\usepackage{graphicx}
\usepackage{bm,amsmath,amssymb,mathrsfs}
\usepackage{amsfonts}
\usepackage{algorithm,algpseudocode,float}
\usepackage{algorithmicx}
\usepackage{multirow}
\usepackage{subfigure}
\usepackage{color}



\makeatletter
\let\OldStatex\Statex
\renewcommand{\Statex}[1][3]{%
  \setlength\@tempdima{\algorithmicindent}%
  \OldStatex\hskip\dimexpr#1\@tempdima\relax}
\makeatother

\def\p{\partial}
\def\DD{\displaystyle}

\def\a{\alpha}
\def\lx{l_x}
\def\ly{l_y}
\def\lpml{l_{\text{pml}}}

\def\N{\mathcal{N}}

\def\star{*}

\title{An Overlapping Domain Decomposition Preconditioner 
  for the Helmholtz equation}

\author{Wei Leng \footnote{State Key Laboratory of Scientific and Engineering Computing,
Chinese Academy of Sciences, Beijing 100190, China. Email: {wleng@lsec.cc.ac.cn}.}
,
Lili Ju \footnote{Department of Mathematics, University of South Carolina,
Columbia, SC 29208, USA. Email: {ju@math.sc.edu}.}
}

\begin{document}

\maketitle

\abstract{In this paper, 
based on the overlapping domain decomposition method (DDM) proposed in \cite{Leng2015},
an one step preconditioner is proposed to solve 2D high frequency Helmholtz equation.  
 %
 The computation domain is decomposed in both $x$ and $y$ directions,
 and the local solution on each subdomain is updated simultaneously in one iteration,
 thus there is no sweeping along certain directions.
 In these ways, the overlapping DDM is similar to the popular DDM for Poisson problem.
 The one step preconditioner simply  take the restricted source on each subdomain,
 solve the local problems and summarize the local solutions on all subdomains including 
 their PML area. 
 The complexity of solving the problem with the preconditioner is $O(N n_{\text{iter}})$,
 where $n_{\text{iter}}$ is the number of iteration, 
 and it is shown numerically that $n_{\text{iter}}$ is 
 proportional to the number of subdomains in one direction.
 %
 2D Helmholtz problem with nearly a billion unknowns 
   are solved efficiently with the preconditioner on massively parallel
   machines.
}




\section{Introduction}

We consider in this paper to solve the Helmholtz equation in the full space $\mathbb{R}^2$, with Sommerfeld radiation condition,
\begin{align} \label{eq:helm}
  \Delta u + k^2 u &= f \qquad  \mbox{in} \,\,\, \mathbb{R}^2, \\
  r^{1/2} (\frac{\p u}{\p r} - \mathbf{i} k u) &\rightarrow 0 \qquad \mbox{as} \,\,\, r = |x| \rightarrow \infty \nonumber 
\end{align}
where $k$ is the wave number.

The domain decomposition method for the Helmholtz equation has been studied for years,
many different DDMs have been proposed based on different boundary conditions 
at the subdomain interface.

The DDM for the Helmholtz equation is very natural. Truncated with perfect match layer,  
the local problem on one subdomain could be approximately solved, 
and the local solution is passed to neighbour subdomains via interface
 to carry on the wave propagation process.
Engquist and Ying \cite{Engquist2011a, Engquist2011b} proposed the sweeping preconditioner
by approximating the inverse of Schur complements in the LDL$^t$ factorization,
and the method has been further developed by Liu and Ying \cite{Liu2015a, Liu2015b} 
using wave addition and dimension recursion. 
Wave traveling in varying medium generates reflections and refractions,   
and it's more reasonable that the DDM  admits reflection at the interface. 
While the sweeping preconditioner \cite{Engquist2011a, Engquist2011b}
use a Dirichlet type interface condition that does not admit reflection, 
a few new DDM is proposed with a reflective interface condition,
such as the source transfer domain decomposition method by Chen and Xiang \cite{Chen2013a,Chen2013b},
Stolk's DDM \cite{Stolk2013}, double sweep preconditioner by Vion an Geuzaine \cite{Vion2014},
and polarized trace method by Zepeda \cite{Zepeda2014}. 
The source transfer DDM \cite{Chen2013a,Chen2013b} admits reflection only in one direction,
recently Du and Wu \cite{Wu2015} modified the method so that it admits reflection 
on both directions. 
Interestingly, we found that the source transfer DDM relates closely to
Stolk's DDM \cite{Stolk2013} and polarized trace method \cite{Zepeda2014},
in the way that choosing the smoothing function in source transfer DDM
to be Heaviside function would lead to an interface condition that 
is similar to the ones of Stolk's DDM and polarized trace method.

The aforementioned domain decomposition methods in the literature usually partition 
the domain into slices in one direction and sweep from one side to the other, 
then sweep backwards. 
Two directions sweeping happens in a recursive way as in Liu and Ying \cite{Liu2015b},
Du and Wu \cite{Wu2015}.
The serial sweeping order causes difficult in scalability in parallel computing, 
and it's impractical to cut too many slices in one direction.

An overlapping DDM is proposed in \cite{Leng2015}.
The popular DDM for problems other than frequency wave domain problem works in such a way that,
the domain is decomposed in multiple directions, and in one step of iteration, 
each subdomain takes the information from its neighbour subdomain, 
update its own solution and prepare the information to be use by
its neighbour subdomains in the next step. 
The overlapping DDM works in the similar way, 
and there is no sweeping along certain direction at all, thus it's suitable for 
large scale parallel computing.
Since the overlapping DDM uses the
source transfer type technique, the reflections is admitted near the interface.
In \cite{Leng2015}, for three layered medium, the reason  is explained why 
the total wave solution is the summation of all incident, reflected and refracted
waves on all subdomains, and the convergence of the method is estimated.
In this paper, we reorganize the overlapping DDM method in a more concise way, and 
an one step preconditioner is proposed. 
Numerical examples are presented to show the preconditioner is simple, 
effective and suitable 
for parallel computing of high frequency wave problems.

The numerical result of the preconditioner shows that the time cost is smaller 
without extra overlapping region. 
However, we still call the preconditioner \textit{overlapping},  
since the PML layer of one subdomain overlap with its neighbours, 
and the solution in the PML layer is added to the total solution,
which is a major difference between this method and the popular DDM for Poisson problem.

The rest of the paper is organized as follows.
In section 2, the overlapping  DDM for Helmholtz equation is
reorganized in a concise way,
and the one step iteration preconditioner is proposed.
In section 3, numerical examples for constant medium, simple layered medium and 
Marmousi model is presented, and the performance of the preconditioner is discussed.

\section{Overlapping domain decomposition preconditioner}

The frequence domain wave equations defined on unbounded domain could
be solved on truncated domain with the perfect matched layer as the absorbing
boundary condition \cite{Berenger, Chew1994}.
To solve Helmholtz problem \eqref{eq:helm}, the unbounded domain $\mathbb{R}^2$
is truncated to a rectangle domain $\Omega = [-l_x-\lpml, l_x+\lpml] \times [-l_y-\lpml, l_y+\lpml]$,
where $\lpml$ is the length of PML layer.
The uniaxial PML method \cite{Chew1994} is used in this paper,
where the complex coordinate is stretched
in $x$ and $y$ direction separately,
$\DD \tilde{x}_j(x_j) = \int_0^{x_j} \sigma_j(t) dt$, $j = 1,2$, and
the PML medium property is chosen that $\sigma_j(t) = 0$ for $ |t| \leq l_j$, and $\sigma_j(t) > 0$
in PML layer $|t| > l_j$. Then the PML equation on the truncated domain is
\begin{equation} \label{eq:pml}
  J^{-1} \nabla \cdot (A \nabla u) + k^2 u = f, \qquad \mbox{in} \,\,\, \Omega,
\end{equation}
where $\DD A(x) = \mbox{diag}\left(\frac{\a_2(x_2)}{\a_1(x_1)}, \frac{\a_1(x_1)}{\a_2(x_2)}\right)$, and $J(x) = \a_1(x_1) \a_2(x_2)$.
The operator of the truncated problem \eqref{eq:pml} is denoted  $\mathcal{L}$,
\begin{equation}
\mathcal{L}(u) =  J^{-1} \nabla \cdot (A \nabla u) + k^2 u.
\end{equation}

The total computation domain is $\Omega$ has 
an interior region $[-\lx, \lx] \times [-\ly, \ly]$, which is parted into 
$Nb_x \times Nb_y$ non-overlapping subdomains. 
Denote $\Delta l = 2 \lx / Nb_x$,  
$x_i = -\lx + i \Delta l$, $i = 0, \ldots, Nb_x$,
and $y_j = -\ly + j \Delta l$, $j = 0, \ldots, Nb_y$,
then the non-overlapping subdomains are $\widetilde{\Omega}_{i,j}
 := [x_i, x_{i+1}] \times [y_j, y_{j+1}]$, 
 $i = 0, \ldots, Nb_x$, $j = 0, \ldots, Nb_y$.

\begin{figure}[ht!]

\subfigure[]{\includegraphics[width=.6\textwidth]{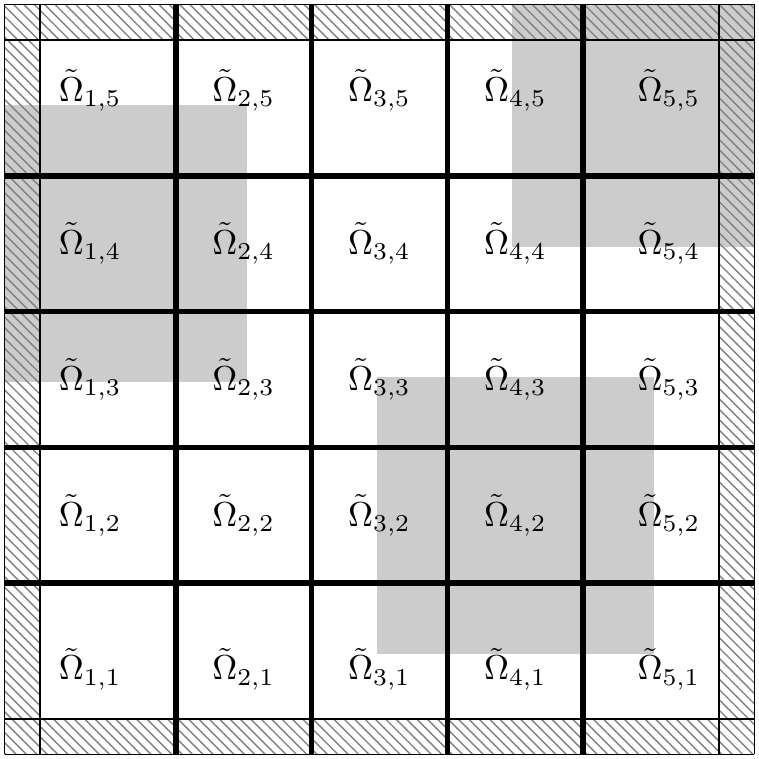} } \quad
\subfigure[]{
             \begin{minipage}[b]{0.2\textwidth}
\includegraphics[scale=1.05]{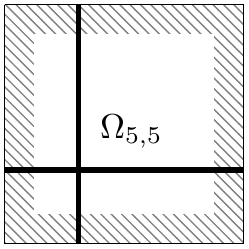} \\
\includegraphics[scale=1.05]{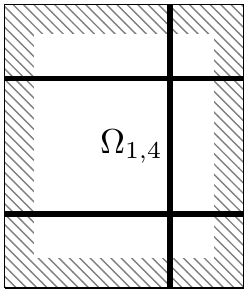} \\
 \includegraphics[scale=1.05]{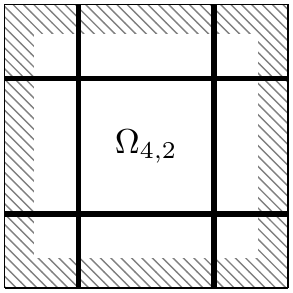}

             \end{minipage} }

\caption{Domain decomposition of 5$\times$5 subdomains. 
The thick lines separates the non-overlapping subdomains $\widetilde{\Omega}_{i,j}$, 
$i,j = 1 \ldots 5$. 
The shadowed areas in (a) are the overlapping subdomain $\Omega_{1,4}$, $\Omega_{5,5}$ and $\Omega_{4,2}$, 
which are also shown in (b).
}
\label{fig:ddm5}
\end{figure}

Then each non-overlapping subdomain $\widetilde{\Omega}_{i,j}$ is extended to overlapping subdomain
$\Omega_{i,j} := 
 [x_{i,0}, x_{i,1}] \times [y_{j,0}, y_{j,1}]$,
where 
\begin{align}
x_{i,0} &= \left\{
\begin{array}{ll}
 x_i - \lpml,     &\,\,\,\,\,\, \quad \mbox{if ~ $i = 0$}\\
 x_i - lo - \lpml,&\,\,\,\,\,\, \quad \mbox{if ~ $i > 0$}\\
\end{array}
\right. \nonumber\\
x_{i,1} &= \left\{
\begin{array}{ll}
 x_{i+1} + lo + \lpml,& \quad \mbox{if ~ $i < Nb_x-1$}\\
 x_{i+1} + \lpml,     & \quad \mbox{if ~ $i = Nb_x-1$}\\
\end{array}
\right.\nonumber\\
y_{j,0} &= \left\{
\begin{array}{ll}
 y_j - \lpml,     &\,\,\,\,\,\, \quad \mbox{if ~ $j = 0$}\\
 y_j - lo - \lpml,&\,\,\,\,\,\, \quad \mbox{if ~ $j > 0$}\\
\end{array}
\right.\nonumber\\
y_{j,1} &= \left\{
\begin{array}{ll}
 y_{j+1} + lo + \lpml,& \quad \mbox{if ~ $j < Nb_y-1$}\\
 y_{j+1} + \lpml,     & \quad \mbox{if ~ $j = Nb_y-1$}\\
\end{array}
\right.\nonumber
\end{align}
and $lo$ is the length of overlapping region.
The domain decomposition with 5$\times$5 subdomains is demonstrated in Fig \ref{fig:ddm5}.

On each subdomain $\Omega_{i,j}$, an local problem with PML layer is set up that solves 
wave field $u_{i,j}$ with given source $f_{i,j}$,
\begin{equation} \label{eq:pmlij}
   J_{i,j}^{-1} \nabla \cdot (A_{i,j} \nabla u_{i,j}) + k^2 u_{i,j} = f_{i,j}, \qquad \mbox{in} \,\,\, \Omega,
\end{equation}
where $J_{i,j}^{-1}$ and $A_{i,j}$ is determined by the PML layer of $\Omega_{i,j}$.
Denote the index set of neighbour $\mathcal{N}_{i,j} = \big\{ (i',j')  \, \big| \, i' \in \{i-1,i,i+1\},  
j' \in \{j-1,j,j+1\},  \,\, \text{and} \,\, (i',j') \neq (i,j)  \big\}$, so the subdomain $\Omega_{i,j}$
has neighbour subdomains $\Omega_{i',j'}$, $(i',j') \in \mathcal{N}_{i,j} $.
Now the overlapping DDM is stated as follows.

On the first step, solve the subdomain problem on $\Omega_{i,j}$ with source restricted 
to interior region $\widetilde{\Omega}_{i,j}$, and the solution is denoted $u_{i,j}^0$, 
\begin{equation} 
J_{i,j}^{-1} \nabla \cdot (A_{i,j} \nabla u_{i,j}^{0}) + k^2 u_{i,j}^{0} 
= f \big|_{\widetilde{\Omega}_{i,j}},
 \qquad \mbox{in} \,\,\, \Omega_{i,j} \label{eq:ddm0} 
\end{equation}

On the successive steps, denote the subdomain solution of step $s$ as $u_{i,j}^s$.
In each step, solve the subdomain problem on $\Omega_{i,j}$ with the residual
of the neighbour subdomains restricted to interior region $\widetilde{\Omega}_{i,j}$ as source,
\begin{equation} 
J_{i,j}^{-1} \nabla \cdot (A_{i,j} \nabla u_{i,j}^{s+1}) + k^2 u_{i,j}^{s+1} 
= - \sum_{(i',j') \in \mathcal{N}_{i,j}} \mathcal{L}(u_{i', j'}^s) \big|_{\widetilde{\Omega}_{i,j}}.
 \qquad \mbox{in} \,\,\, \Omega_{i,j} \label{eq:ddms} 
\end{equation}
Such iteration goes on until the residual is small enough. And the DDM solution is 
\begin{equation}  
u_{\text{DDM}} = \sum_{s>=0} \sum_{i,j} u_{i,j}^s \label{eq:sum}
\end{equation}

The main idea of the overlapping DDM is as follows. On each subdomain $\Omega_{i,j}$, 
the source, e.g. $f_{i,j}$ on $\widetilde{\Omega}_{i,j}$ cause a local wave field $u_{i,j}$, and the residual 
$r_{i,j} := f_{i,j} - \mathcal{L} u_{i,j}$ satisfied that 
$r_{i,j} = 0$ in $\widetilde{\Omega}_{i,j}$, 
and $r_{i,j} \neq 0$ in PML layer of $\Omega_{i,j}$,
thus the residual in PML layer contains the wave field information 
that can be used as incident wave field for neighbour subdomains
to carry on the wave propagation precess.
All subdomains solve the local problem in parallel, 
and send the wave information to neighbour subdomains
in one iteration. After $m$ iterations, the wave have 
approximately propagated over $m$ subdomains.
If there are medium discontinuities in the subdomains, 
reflections will be passed back by the residual in the PML layer 
in the next iteration.
Such wave propagation precess goes on during the iteration, 
and the summation of all incident, reflected and refracted waves
on all subdomains is the total solution. 
Detailed discussion for two subdomains 
with three layered medium could be found in \cite{Leng2015}.

\begin{figure}[ht!]
\centering
 
\subfigure[]{\includegraphics[scale=1.1]{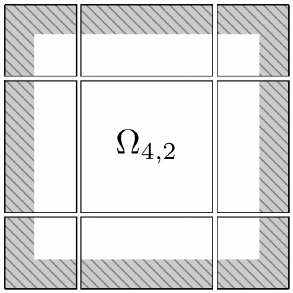}   }  \quad
\subfigure[]{\includegraphics[scale=1.1]{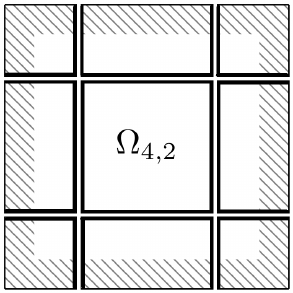}  } \quad
\subfigure[]{\includegraphics[scale=1.1]{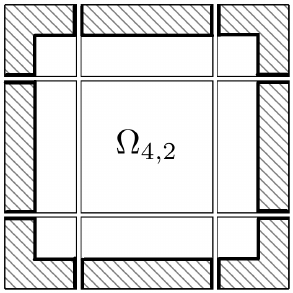}  }

\caption{Illustration of the information passing from $\Omega_{4,2}$ to its neighbors.
The subdomain $\Omega_{4,2}$ is separated into nine parts that belong to different 
non-overlapping subdomains, 
namely $\widetilde{\Omega}_{4,2}$ and $\Omega_{4,2} \big|_{\widetilde{\Omega}_{i',j'}}$, where $(i',j') \in \N_{4,2}$.
(a) The information is contained in the residual in $\Omega_{4,2}$'s PML region.  
(b) The information is contained in the incident boundaries, marked with thick lines.
(c) The information is contained in the incident boundaries, also marked with thick lines.
}
\label{fig:info}
\end{figure}

%
%
%
%
%
%
%
To explain why the overlapping DDM works well for domain decomposition in both $x$ and $y$  
directions,
we elaborate on the wave information passing 
 from $\Omega_{i,j}$ to its neighbor  $\Omega_{i',j'}$, 
$(i',j') \in \N(i,j)$. By \eqref{eq:ddms}, the wave information
is passed with the residual in the subdomain's PML layer, as shown in 
Fig \ref{fig:info}-(a). Alternately, the wave field in the 
subdomain's PML layer could be recovered using only the
values on incident boundaries, as shown in Fig \ref{fig:info}-(b) and (c), 
and the residual in the PML layer is then recovered.
In either case, the wave information is passed 
not only in $x$ direction and  $y$ direction, but also 
in the corner direction.
The amount of the information passed to corner neighbour subdomais increases
as the length of PML layer or the length of extra overlap region increases.

%
%
%
%
%
%
%
The overlapping DDM is more effective when using as a preconditioner than a solver.
The one step preconditioner is chosen, which simply 
 solves the subdomain problem on $\Omega_{i,j}$ with source restricted 
to interior region $\widetilde{\Omega}_{i,j}$, 
\begin{equation} 
J_{i,j}^{-1} \nabla \cdot (A_{i,j} \nabla u_{i,j}^{0}) + k^2 u_{i,j}^{0} 
= f \big|_{\widetilde{\Omega}_{i,j}}.
 \qquad \mbox{in} \,\,\, \Omega_{i,j} \label{eq:ddm0} 
\end{equation}
and summarizes the solutions on all subdomains to get the approximate solution 
used in the preconditioning procedure , 
\begin{equation} \label{eq:pc}
u_{\text{PC}} =  \sum_{i,j} u_{i,j}^0
\end{equation}

\textbf{Remark:}
A smoothing function could be multiplied to the solution in the PML region 
to keep the solution in $H^1(\Omega)$, such technique is used in \cite{Chen2013a}. 
In \eqref{eq:ddms}, 
$\mathcal{L}(u_{i',j'}^s)$ could be substituted with $\mathcal{L}(\beta_{i',j'} u_{i',j'}^s)$,
where $\beta_{i,j}$ is a smoothing function for subdomain $\Omega_{i,j}$, 
such that $\beta_{i,j} \in C^2(\mathbb{R}^2)$, 
$\beta_{i',j'} \big|_{\widetilde{\Omega}_{i',j'}} = 1$,
and $\beta_{i',j'} \big|_{R^2 \char`\\ {\Omega}_{i',j'}} = 0$.
Mean while in \eqref{eq:sum} and \eqref{eq:pc}, the summation over $u^s_{i,j}$ 
could be   substituted with $\beta_{i,j} u^s_{i,j}$.
In our numerical experiments, since the remaining value near the  PML layer outside boundary
is negligible, we simply omit the smoothing function.

\section{Numerical experiments}

Three numerical experiments that includes constant medium, simple layered medium and
Marmousi model are carried out to test the performance of the overlapping DDM.

Finite difference method with second order  accuracy is used to
discretize the Helmholtz equation. 
The PML layer is of 30 grid points width by default.
Single shot in the subdomain $\Omega_{0,0}$ is taken as the source, 
and the position is $(x^s, y^s)$ where $\DD x^s = x_0 + \frac{1}{4} \, \Delta l$, 
and $\DD y^s = y_0 + \frac{1}{3} \, \Delta l$.     
 The shape of the shot is an approximate delta function,
 $\DD f_{i,j} = \frac{1}{h_x h_y}\delta(i h_x - x^s, j h_y - y^s)$, 
  where $h_x$, $h_y$
 are the grid size in $x$ and $y$ direction, respectfully.
The relative tolerance of linear solver is 10$^{-10}$ for all cases.

The number of equivalent sweeping, defined as $\DD \frac{n_{\text{iter}}}{\max\{Nb_x, Nb_y\}}$, 
is used to measure the effectiveness of the preconditioner.
A source at one side of the domain generates the wave that pass to the other side,
and it takes at least number of subdomains in the direction to accomplish the traveling.
Thus, we are satisfied if  the number of equivalent sweeping does not change 
if number of subdomains grows.

The Tianhe-2 cluster is used in our numerical experiments, each node of 
the cluster includes two 2.2GHz Xeon E5-2692 processors with 12 cores. 
The number of processors in use are $Nb_x \times Nb_y$.

\subsection{Constant medium}

The overlapping DDM is tested for constant $k = 1$ on square domain $[0,1]\times[0,1]$. 
First, we fixed the number of subdomains, and increase both the problem size and the wave number,
to see the effect of increasing frequency to the preconditioner. 
Note that the number of PML layer points $n_{\text{pml}}$ increases
as problem size, so that on each subdomain the ratio of the PML area to the total area  is fixed. 
$8 \times 8$ subdomains in $x$ and $y$ direction is used, 
and on each subdomain the local problem is factorized and solved with direct solver.

\begin{table}[ht!]
\centering
\begin{tabular}{|c|c|c|c|c|c|c|c|}
  \hline
  Size       & $n_{\text{pml}} $ & Freq            & No.            & No. equ & Time & Time \\
             &                   & $\omega / 2\pi$ & \,\, Iter \,\, & sweep   & Fact & Iter \\
\hline
  800$^2$    & 30                & 71.6            & 53             & 7       & 0.41 & 2.53 \\
  1,600$^2$  & 60                & 143             & 42             & 5       & 1.39 & 7.55 \\
  3,200$^2$  & 120               & 287             & 37             & 5       & 6.42 & 32.1 \\
  6,400$^2$  & 240               & 573             & 34             & 4       & 33.7 & 111  \\
  12,800$^2$ & 480               & 1147            & 31             & 4       & 192  & 440  \\
  \hline
\end{tabular}
\caption{The performance of the preconditioner for const medium problem, 
with number of subdomains fixed.} \label{tab:const1}
\end{table}

The result is shown in table \ref{tab:const1}. 
The time cost for factorizing and solving local problem on the subdomain increases
as the local problem size grows. 
The number of equivalent sweeping is around 5 for different frequency, thus 
the preconditioner is not affected by the frequency of the problem in the constant medium case.
Actually, the number of equivalent sweeping decrease a little bit as the problem size 
 grows, which is expected since the increasing number of 
 points in PML layer  leads to better absorption at the subdomain boundary.

Second, the weak scalability test is preformed. 
Both the problem size and the number of subdomains increase,
while the subdomain problem size is fixed.
Such test exams whether a method is suitable for large scale parallel computing.   
The result is shown in Table \ref{tab:const2}. 
Since the problem size of each subdomain is fixed, a fixed 
setup time around 2.2s to factorize the local problem is required.  
The number of equivalent sweeping is kept almost unchanged as the 
total problem size increases, however, we found that 
the number of GMRES restart need to grow to maintain such iteration numbers.  
The number of iteration doubles as the problem size in one direction doubles,
cause the total solving time doubles, and the iteration time is far larger 
than the fixed setup time as the problem grows large.

\begin{table}[ht!]
\centering
\begin{tabular}{|c|c|c|c|c|c|c|}
  \hline
 Size       & $Nb_x \times Nb_y$ & Freq            & GMRES   & No.            & No. equ & Time     \\
            &                    & $\omega / 2\pi$ & Restart & \,\, Iter \,\, & sweep   & solve(s) \\
\hline
 600$^2$    & 2 $\times$ 2       & 55              & 30      & 9              & 5       & 3.5      \\
 1,200$^2$  & 4 $\times$ 4       & 105             & 30      & 23             & 6       & 7.3      \\
 2,400$^2$  & 8 $\times$ 8       & 205             & 30      & 52             & 7       & 16       \\
 4,800$^2$  & 16 $\times$ 16     & 405             & 60      & 106            & 7       & 33       \\
 9,600$^2$  & 32 $\times$ 32     & 805             & 120     & 213            & 7       & 79       \\
 19,200$^2$ & 64 $\times$ 64     & 1605            & 240     & 425            & 7       & 324      \\
  \hline
\end{tabular}
\caption{The performance of the preconditioner for const medium problem, 
with subdomain problem size fixed.} \label{tab:const2}
\end{table}

\subsection{Layered medium}
\begin{figure}[ht!]
\centerline{
\includegraphics[width=.5\textwidth]{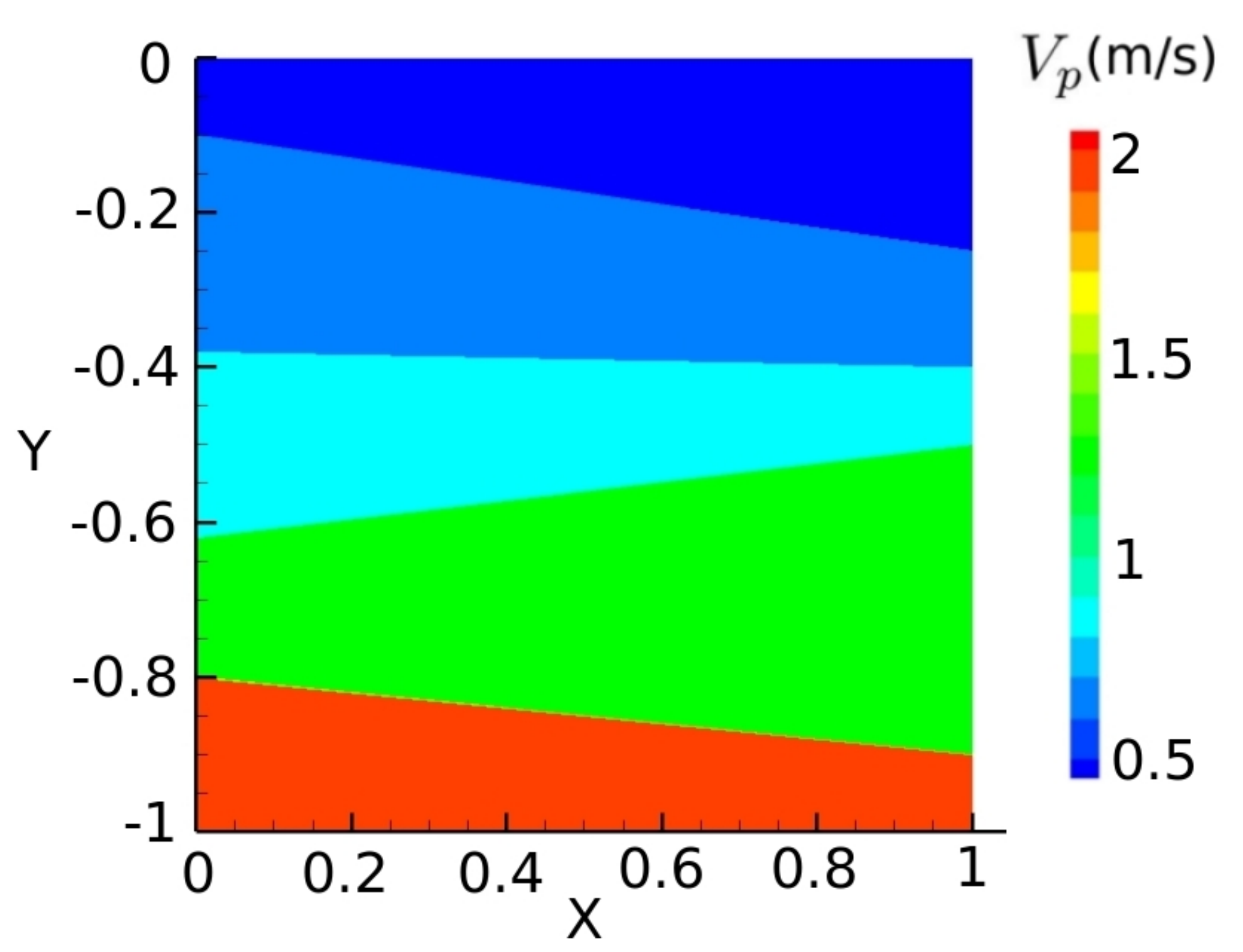} \quad
\includegraphics[width=.5\textwidth]{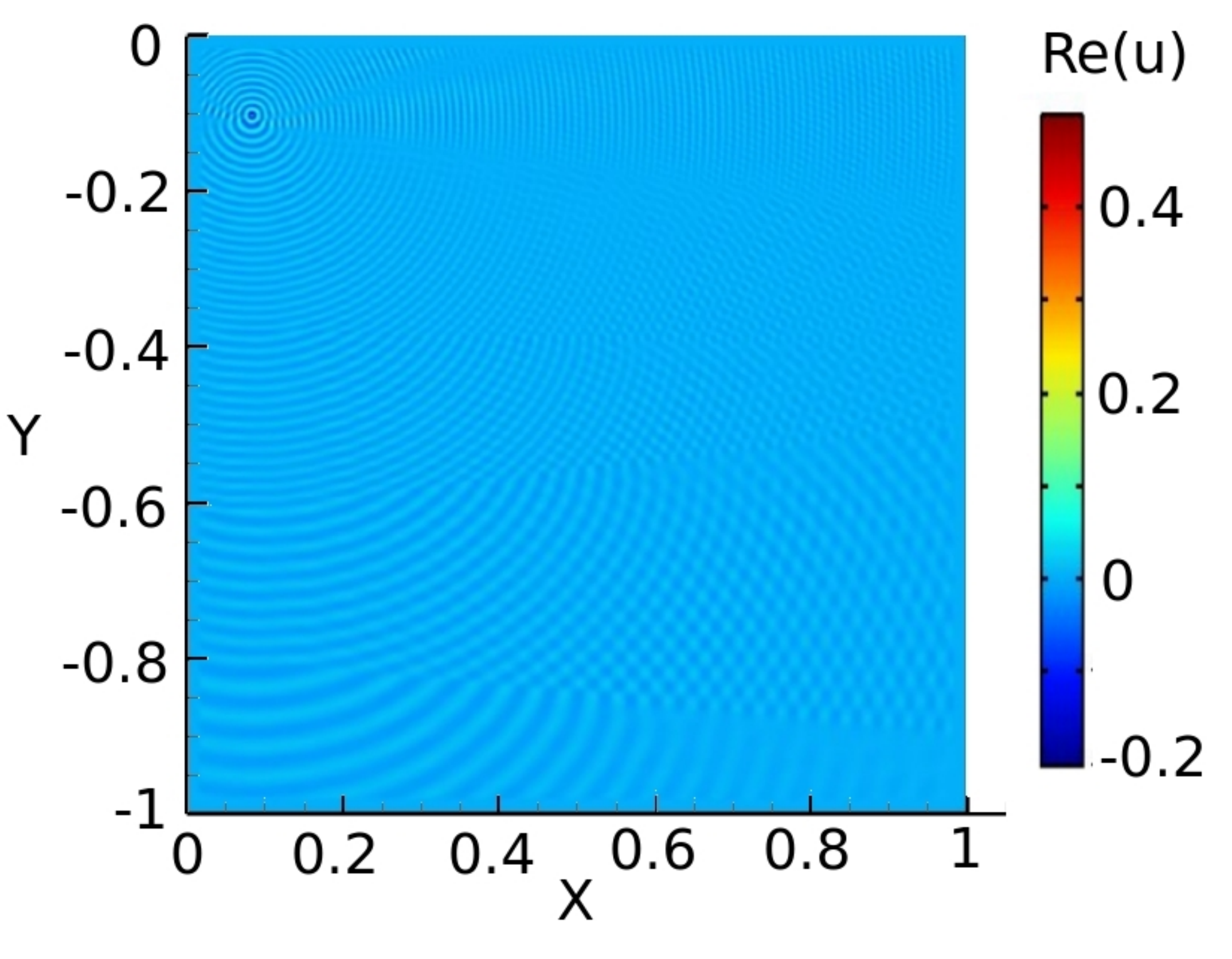} 
}
\caption{Left: layered medium velocity profile. Right: numerical solution of problem size $1200^2$. 
}
\label{fig:5layer}
\end{figure}

The overlapping DDM is tested for simple layered medium problem on square domain $[0,1]\times[-1,0]$,
where five layers of medium lie nearly horizontally, as shown in Fig \ref{fig:5layer} left.  
The numerical solution of the problem with problem size $1200^2$ is shown in Fig \ref{fig:5layer} right.

Again, we fixed the number of subdomains, and increase both the problem size and the wave number,
to see the effect of increasing frequency to the preconditioner. 
$8 \times 8$ subdomains in $x$ and $y$ direction is used, and 
the result is shown in table \ref{tab:layer1}. 
The number of equivalent sweeping is around 10 for different frequency,
a bit larger than that of constant medium problem in similar test,  
thus the preconditioner is not affected by the frequency of the problem in this case. 
Similar to the result of const medium problem, 
the number of equivalent sweeping decrease a little bit as the problem size 
grows.   

\begin{table}[ht!]
\centering
\begin{tabular}{|c|c|c|c|c|c|c|c|}
  \hline
  Size       & $n_{\text{pml}} $ & Freq            & No.            & No. equ & Time & Time \\
             &                   & $\omega / 2\pi$ & \,\, Iter \,\, & sweep   & Fact & Iter \\
\hline
  800$^2$    & 30                & 35.8            & 88             & 11      & 0.49 & 4.70 \\
  1,600$^2$  & 60                & 71.6            & 83             & 10      & 1.79 & 18.3 \\
  3,200$^2$  & 120               & 143             & 78             & 10      & 8.18 & 79.4 \\
  6,400$^2$  & 240               & 287             & 74             & 9       & 34.0 & 231  \\
  12,800$^2$ & 480               & 573             & 70             & 9       & 192  & 935  \\
  \hline
\end{tabular}
\caption{The performance of the preconditioner for simple layered model, 
with number of subdomains fixed. } \label{tab:layer1}
\end{table}

The weak scalability test is also preformed, 
and the result is shown in Table \ref{tab:layer2}.
The number of extra overlapping points $N_{OL} = 0, 50$ is tested 
to evaluate the effectiveness of enlarging overlapping region.
The fixed setup time to factorize the local problem is around 2.2s for $N_{OL} = 0$
and 3.1s for $N_{OL} = 50$, respectfully.
The number of equivalent sweeping is kept almost unchanged as the 
total problem size increase.  
Enlarging the overlapping region results smaller number of iteration,
however since the local problem size increases and it take longer time to solve the local problem, 
the total time cost is 
bigger. So in this case, it's better without extra overlapping region.

\begin{table}[ht!]
\centering
\begin{tabular}{|c|c|c|c|c|c|c|c|}
\hline
 Size                                        & $Nb_x \times Nb_y$                  & Freq                       & $N_{OL}$ & GMRES   & No.            & No. equ & Time     \\
                                             &                                     & $\omega / 2\pi$            &          & Restart & \,\, Iter \,\, & sweep   & solve(s) \\
\hline
 \multirow{2}{\star}{600$^2$}       & \multirow{2}{\star}{2 $\times$ 2}   & \multirow{2}{\star}{27.5}  & 0        & 30      & 19             & 10      & 4.86     \\
                                             &                                     &                            & 50       & 30      & 16             & 8       & 8.28     \\
\hline
 \multirow{2}{\star}{1,200$^2$}   & \multirow{2}{\star}{4 $\times$ 4}   & \multirow{2}{\star}{52.5}  & 0        & 30      & 46             & 12      & 12.4     \\
                                             &                                     &                            & 50       & 30      & 41             & 10      & 22.8     \\
\hline
 \multirow{2}{\star}{2,400$^2$ }   & \multirow{2}{\star}{8 $\times$ 8}   & \multirow{2}{\star}{102.5} & 0        & 30      & 88             & 11      & 24.2     \\
                                             &                                     &                            & 50       & 30      & 77             & 10      & 41.2     \\
\hline
 \multirow{2}{\star}{4,800$^2$ }   & \multirow{2}{\star}{16 $\times$ 16} & \multirow{2}{\star}{202.5} & 0        & 60      & 173            & 11      & 53.3     \\
                                             &                                     &                            & 50       & 60      & 153            & 10      & 91.9     \\
\hline
 \multirow{2}{\star}{9,600$^2$ }   & \multirow{2}{\star}{32 $\times$ 32} & \multirow{2}{\star}{402.5} & 0        & 120     & 338            & 11      & 127      \\
                                             &                                     &                            & 50       & 120     & 294            & 9       & 198      \\
\hline
 \multirow{2}{\star}{19,200$^2$} & \multirow{2}{\star}{64 $\times$ 64} & \multirow{2}{\star}{802.5} & 0        & 240     & 680            & 11      & 560      \\
                                             &                                     &                            & 50       & 240     &  615              & 10      &  855     \\
\hline
\end{tabular}
\caption{The performance of the preconditioner for simple layered model, 
with subdomain problem size fixed.} \label{tab:layer2}
\end{table}

\subsection{Marmousi model}

\begin{figure}[ht!]
\begin{center}
\includegraphics[width=.9\textwidth]{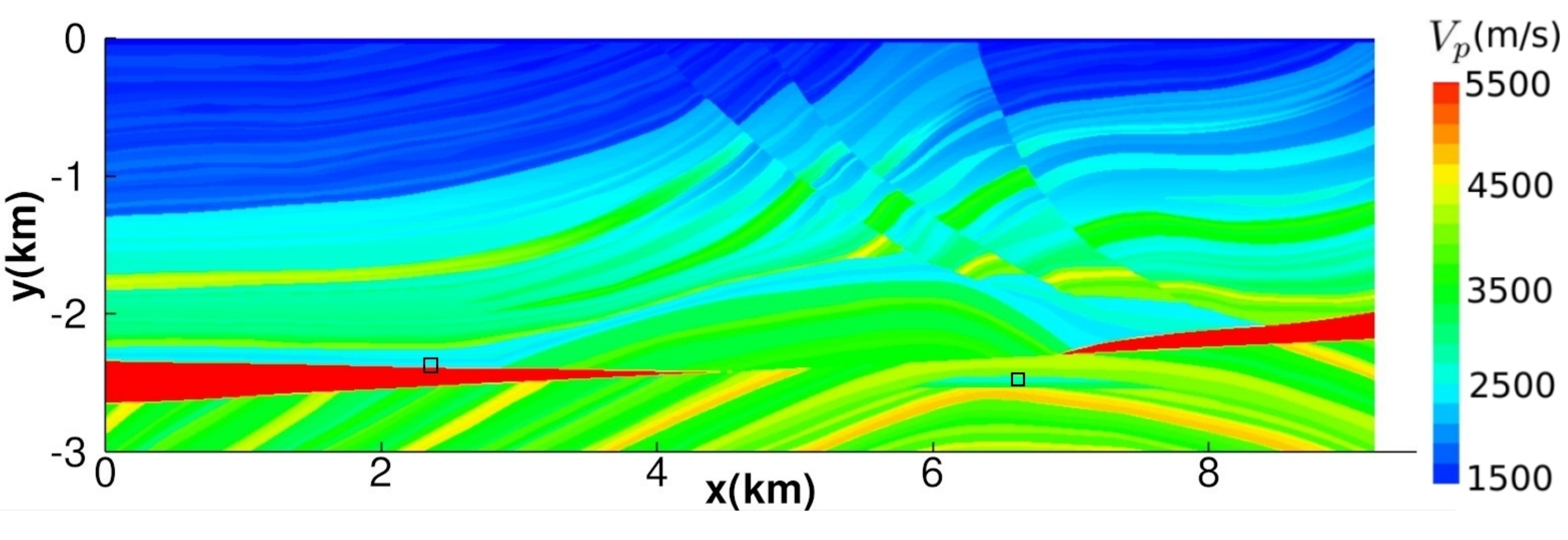}
\end{center}
\caption{Velocity profile of Marmousi model.}
\label{fig:mar-vel}
\end{figure}

\begin{figure}[ht!]
\begin{center}
\includegraphics[width=.9\textwidth]{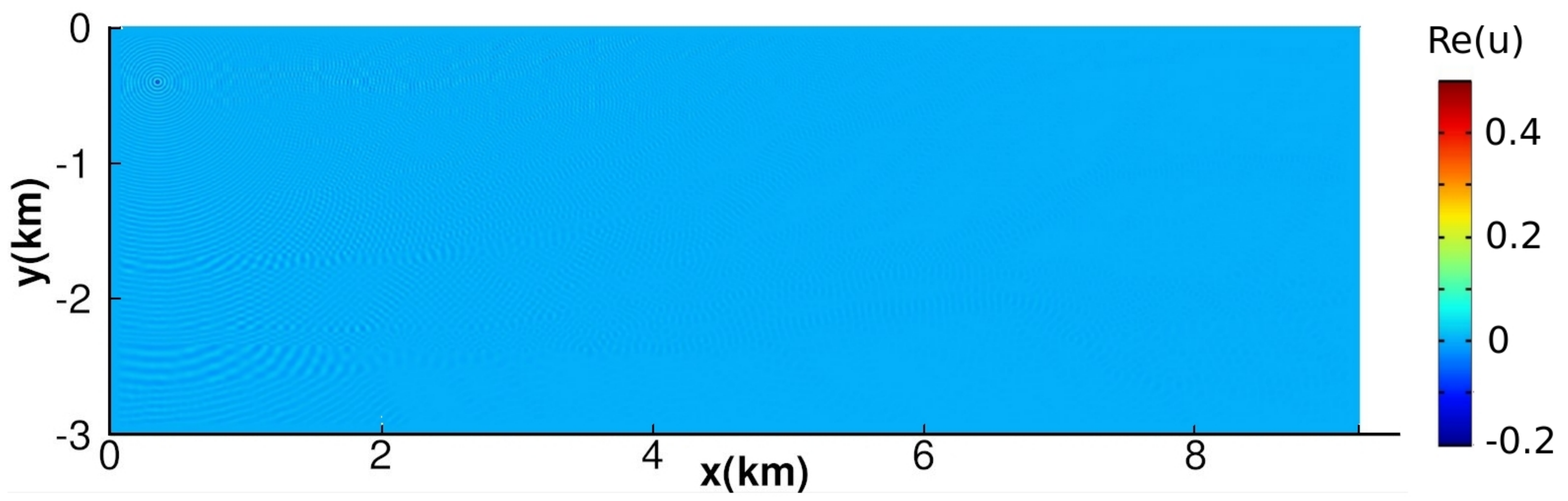}
\end{center}
\caption{Solution of problem size 4,275 $\times$ 1,425.}
\label{fig:mar-sol}
\end{figure}

At last, the preconditioner is tested on the 2D Marmousi model in seismology,
which is  $3,000$ m deep and $9,200$ m wide.  Only P-wave
is considered, thus elastic wave equation becomes an acoustic equation.
The velocity profile is shown in Fig \ref{fig:mar-vel}, the
maximum velocity is 5500 km/s and the minmum velocity is 1500 km/s.
The numerical solution of the problem with problem size $4,275\times 1,425$ is
 shown in Fig \ref{fig:mar-sol} right.

The weak scalability test is also preformed, 
and the result is shown in Table \ref{tab:marm}.
The number of equivalent sweeping is kept almost unchanged as the 
total problem size increases. 
The largest size problem has 977,407,500 unknowns, and 4,332 processors is used to solve it.

\begin{table}[ht!]
\centering
\begin{tabular}{|c|c|c|c|c|c|c|}
  \hline
 Size                    &  $Nb_x \times Nb_y$  &             Freq  &  Restart  &             No.  &  No. equ  &   Time  \\
                         &                      &  $\omega / 2\pi$  &           &  \,\, Iter \,\,  &    sweep  &  solve  \\
\hline
 4,275 $\times$ 1,425    &  9 $\times$ 3        &             63.5  &       30  &              86  &       10  &   56.6  \\
 8,550 $\times$ 2,850    &  18 $\times$ 6       &              123  &       30  &             151  &        8  &    111  \\
 17,100 $\times$ 5,700   &  36 $\times$ 12      &              242  &       30  &             290  &        8  &    230  \\
 34,200 $\times$ 11,400  &  72 $\times$ 24      &              479  &       60  &             697  &       10  &    605  \\
 54,150 $\times$ 18,050  &  114 $\times$ 38     &              756  &      120  &           1,123  &       10  &   1177  \\
\hline
\end{tabular}
\caption{The performance of the preconditioner for Marmousi model, 
with subdomain problem size fixed.} \label{tab:marm}
\end{table}

%
%

\section*{Acknowledgments}
This work is supported by the National 863 Project of China
under the grant number 2012AA01A309, and the National Center for
Mathematics and Interdisciplinary Sciences of the Chinese Academy of
Sciences.


\end{document}